\documentclass[%
twocolumn
]{mpi2015-cscpreprint}

\usepackage{graphicx}      

\usepackage[american]{babel}

\usepackage{todonotes}

\usepackage{geometry}
\usepackage{graphicx} 
\usepackage{float}
\usepackage{color}
\usepackage{url}
\usepackage{amsmath}
\usepackage{amsfonts}
\usepackage{amssymb}
\usepackage{scalefnt}
\usepackage{bbm}

\usepackage[linesnumbered, ruled, vlined]{algorithm2e}

\newcommand{\bA}{{\mathbf A}}
\newcommand{\bB}{{\mathbf B}}
\newcommand{\bC}{{\mathbf C}}
\newcommand{\bD}{{\mathbf D}}
\newcommand{\bE}{{\mathbf E}}

\newcommand{\bI}{{\mathbf I}}
\newcommand{\bM}{{\mathbf M}}

\newcommand{\bX}{{\mathbf X}}
\newcommand{\bY}{{\mathbf Y}}

\newcommand{\bx}{{\mathbf x}}

\def\IR{{\mathbb R}}
\def\IC{{\mathbb C}}
\newcommand{\cD}{ {\cal D} }

\def\IL{{\mathbb L}}
\def\IV{{\mathbb V}}
\def\IW{{\mathbb W}}

\newcommand{\IM}{{{{\mathbb L}_s}}}

\newcommand{\bLambda}{\boldsymbol{\Lambda}}

\newcommand{\bhA}{{\hat{\textbf A}}}
\newcommand{\bhB}{{\hat{\textbf B}}}
\newcommand{\bhC}{{\hat{\textbf C}}}

\newcommand{\bhE}{{\hat{\textbf E}}}

\newcommand{\bhx}{{\hat{\textbf x}}}

\newcommand{\hH}{{\hat{H}}}

\newtheorem{lemma}{Lemma}

\newcommand{\bone}{{\mathbbm{1}}}

\begin{document}
	
	\title{Approximating a flexible beam model in the Loewner framework}

	\author[$\ast$]{Alexander Zuyev}
\affil[$\ast$]{Otto von Guericke University Magdeburg, Germany, Max Planck Institute for Dynamics of Complex Technical Systems, Magdeburg, Germany, and the Institute of Applied Mathematics \& Mechanics, National Academy of Sciences of Ukraine.\authorcr
	\email{zuyev@mpi-magdeburg.mpg.de}, \orcid{0000-0002-7610-5621}}
	
	\author[$\dagger$]{Ion Victor Gosea}
	\affil[$\dagger$]{Max Planck Institute for Dynamics of Complex Technical Systems, Magdeburg, Germany.\authorcr
		\email{gosea@mpi-magdeburg.mpg.de}, \orcid{0000-0003-3580-4116}}

	\shorttitle{Approximating a flexible beam model in the Loewner framework}
	\shortauthor{A. Zuyev and I. V. Gosea}
	\shortdate{}
	

\abstract{
The paper develops the Loewner approach for data-based modeling of a linear distributed-parameter system. This approach is applied to a controlled flexible beam model coupled with a spring-mass system. The original dynamical system is described by the Euler-Bernoulli partial differential equation with the interface conditions due to the oscillations of the lumped part. The transfer function of this model is computed analytically, and its sampled values are then used for the data-driven design of a reduced model. A family of approximate realizations of the corresponding input-output map is constructed within the Loewner framework. It is shown that the proposed finite-dimensional approximations are able to capture the key properties of the original dynamics over a given range of observed frequencies. The robustness of the method to noisy data is also investigated.}

\maketitle


\section{Introduction}

Model reduction techniques can be employed to replace a large-scale system with a complex structure (characterized by multidimensional systems of ordinary differential equations and/or partial differential equations), with a much simpler and smaller dynamical system (characterized by few equations with well-understood dynamics). In the last decades, there have been many methodologies proposed in this direction; we refer the reader to \cite{ACA05,morQuaR14,kutz2016dynamic,BOCW17,ABG20} for more details.

A viable alternative to using classical model reduction approaches based on single or double-sided projections (that usually require explicit access to a large-scale model) is to use instead data-driven methods. These latter do not require explicit access to the large-scale model's structure or matrices. We mention here the Loewner framework (LF) \cite{MA07},  Vector fitting (VF) \cite{morGusS99}, or the AAA algorithm \cite{NST18}. When using these, low-order models can be constructed directly from data in the frequency domain (samples of the transfer function). Such methods can be viewed as rational approximation tools by means of interpolation (LF), least-squares fit (VF), or a mixed approach (AAA). Other data-driven methods that has emerged in recent years are dynamic mode decomposition (DMD) and operator inference (OpInf), which use time-domain snapshots of the state variables and then fit a particular structured model by computing the appropriate matrices (in reduced coordinates). Details on DMD can be found in \cite{kutz2016dynamic}, while details on OpInf can be found in \cite{morPehW16,morBenGKetal20}.

We consider here the problem of data-driven rational approximation by means of fitting a linear time-invariant (LTI) dynamical system to a set of measurements (in the frequency domain). The fitted LTI system is characterized in the state-space by the following equations:
\begin{equation}\label{eq:linsys}
\begin{cases}
\bE \dot{\bx}(t)=\bA\bx(t)+\bB u(t),\\ y(t) \ \hspace{+1.5mm} =\bC\bx(t)+\bD u(t),
\end{cases}
\end{equation}
where $u(t)\in \mathbb R$ is the input, $y(t)\in\mathbb R$ is the output, $\bx(t)\in \mathbb R^n$ is the state vector, and the system matrices are  $\bA, \bE \in\IR^{n\times n},~\bB,\bC^T\in\IR^{n\times 1}$. The transfer function of (\ref{eq:linsys}) is given by $H(s) = \bC(s\bE-\bA)^{-1}\bB+\bD$. We refer to \cite{ACA05} for more details on various methodologies especially tailored to the reduction of linear systems. In recent years, such methods for linear systems have been steadily extended to particular classes of structured linear systems \cite{beattie2009interpolatory,duff2015realization}, or even to nonlinear structured systems in \cite{gosea2022iterative} (without preservation of structure) or in \cite{benner2021structure} (with preservation of structure). The structures treated include distributed parameters, delay terms or integro-differential equations. In a data-driven setup, the structure-preserving approach is \cite{schulze2016data} was proposed. Note that the analytical representations of transfer functions have been obtained only for particular classes of distributed parameter systems. We refer to \cite{curtain2009transfer,altiner2019modeling} for surveys of results in this area. In the former tutorial article, the authors provide the derivation of a variety of (irrational) transfer functions for systems described by partial-differential equations. It is also shown that the choice of boundary conditions have an influence on the dynamics and on the locations of poles and zeros. In most practical situations, it is desirable to approximate the irrational transfer function  by a rational one, for the purpose of controller design.

The Loewner framework approach was shown to be extremely powerful in data-based control problems for wide classes of finite-dimensional control systems, whose transfer functions are rational \cite{ALI17,ABG20}.
However, the efficiency of the Loewner framework for distributed-parameter systems (characterized by irrational transfer functions) still remains to be verified, and the present paper aims at filling this gap. Preliminary analysis was provided for linear time-delay systems in \cite{schulze2018data,LM18}, fractional-order systems in \cite{casagrande2019integer}, or for control purposes in \cite{gosea2021loewner,poussot2022interpolation}. A recent overview was provided in \cite{karachalios2021loewner}, including amongst others, rational approximation of the Bessel function, of a hyperbolic sine, and of a vibrating beam model from \cite{curtain2009transfer}. Hence, the LF was studied in the context of approximating infinite-dimensional models of vibrating beams (with finite-dimensional ones). However, as far as the authors are aware, this is the first contribution that also takes into account the effects of perturbed data, i.e., under additive Gaussian noise.


\section{Vibrating beam with attached mass}
Consider the Euler--Bernoulli equation describing the transverse vibrations of a flexible beam of length $l$:
\begin{equation}\label{EulerBernoulli}
\ddot w(x,t)+\frac{EI}{\rho} w''''(x,t)+d \, \dot w''''(x,t)=\frac{1}{\rho}\sum\limits_{j=1}^k \psi_j''(x) u_j,
\end{equation}
where $w(x,t)$ is the beam deflection at point $x\in[0,l]$ and time $t$,
$E$ is the Young's modulus, $I$ is the area moment of inertia of the cross-section, $\rho$ is the mass per unit length of the beam,
and $d\ge 0$ is the structural damping coefficient.
We denote the derivative with respect to time by a dot,
while the prime denotes the spatial derivative (i.e., with respect to $x$).
We assume that a mass-spring system (shaker) is attached to the beam at point $x=l_0$, so that equation~\eqref{EulerBernoulli} holds for $x\in [0,l_0)$ and $x\in (l_0,l]$, and the interface condition is imposed at $x=l_0$:
\begin{equation}\label{interface}
(m\ddot w+\varkappa w)\Big|_{x=l_0}=(EIw'')'\Big|_{x=l_0-0}-(EIw'')'\Big|_{x=l_0+0}+u_0.
\end{equation}
The beam is hinged at both ends, which is formalized by the boundary conditions
\begin{equation}\label{BC}
w\Big|_{x=0}=w\Big|_{x=l}=0,\quad w''\Big|_{x=0}=w''\Big|_{x=l}=0.
\end{equation}
System~\eqref{EulerBernoulli}--\eqref{BC} is controlled by the force $u_0$ applied to the shaker at $x=l_0$ and $k$ piezo actuators, whose actions $u_j$ are characterized in terms of shape functions $\psi_j(x)$, $j=1, ... ,k$.
It is also assumed that $p$ piezo sensors are located at the points $x=l_i$, $i=1,..,p$, i.e. the system outputs are
\begin{equation}
y_i(t) = \left.\frac{\partial^2 w(x,t)}{\partial x^2}\right|_{x=l_i},\quad i=1,...,p.
\label{outputs}
\end{equation}

The above mathematical model has been presented in~\cite{KZB2021,KZ2021} for the case without damping;
here we take into account the structural damping by introducing the parameter $d$ in~\eqref{EulerBernoulli}.

\section{Computation of the transfer function}
\label{sec_tf}
Let $u_0(t)$, $u_1(t)$, ..., $u_k(t)$ $(t\ge 0)$ be inputs of the control system~\eqref{EulerBernoulli}--\eqref{outputs} with zero initial data,
denote the Laplace transform of the inputs and outputs by
$$
U_j(s) = \int_0^{+\infty} u_j(t) e^{-st}dt,\quad j=0,1,...,k,
$$
and
$$
Y_i(s) = \int_0^{+\infty} y_i(t) e^{-st}dt,\quad i=1,...,p,
$$
respectively.
%


After introducing the Laplace transform of $w(x,t)$ with respect to $t$: $W_1^s(x)=\int_0^{+\infty} w(x,t) e^{-st}dt$, we obtain from~\eqref{EulerBernoulli} the following system of ordinary differential equations for $W^s(x)=\left(W_1^s,W_2^s,W_3^s,W_4^s\right)^T$:
\begin{equation}\label{ODE_vector}
\frac d{dx}W^s(x)=A W^s(x) + \Phi^s(x), A=\left(
\begin{array}{cccc}
0 & 1 & 0 & 0 \\
0 & 0 & 1 & 0 \\
0 & 0 & 0 & 1 \\
- 4 \gamma^4 & 0 & 0 & 0 \\
\end{array}
\right),
\end{equation}
with
\begin{equation}
\Phi^s(x)=\left(\begin{array}{c}0 \\ 0 \\ 0 \\ \phi^s(x)\end{array}\right),\; \gamma = \frac{\alpha\sqrt{2s}}{2},
\end{equation}
\begin{equation}
\alpha^4 = \frac{\rho}{EI+\rho d s}>0,\, \phi^s(x) = \frac{1}{EI+\rho d s}\sum_{j=1}^k \psi_j(x) U_j(s),
\label{phi}
\end{equation}
and the $s$ variable is treated as a parameter in~\eqref{ODE_vector}.
The general solution of~\eqref{ODE_vector} is represented with the matrix exponential as follows:
\begin{equation}\label{GenSol}
W^s(x)=\left\{
\begin{array}{ll}
e^{xA}\bar W^0 + \int_0^x e^{(x-y)A}\Phi^s(y)dy, & x\in[0,l_0], \\
e^{(x-l)A}\bar W^l - \int_x^l e^{(x-y)A}\Phi^s(y)dy, & x\in(l_0,l],
\end{array}
\right.
\end{equation}
where
\begin{equation*}
e^{xA}=\left(
\begin{array}{cccc}
z_1(x) & z_2(x) & z_3(x) & z_4(x) \\
-4\gamma^4 z_4(x) & z_1(x) & z_2(x) & z_3(x) \\
-4\gamma^4 z_3(x) & -4\gamma^4 z_4(x) & z_1(x) & z_2(x) \\
-4\gamma^4 z_2(x) & -4\gamma^4 z_3(x) & -4\gamma^4 z_4(x) & z_1(x) \\
\end{array}
\right),
\end{equation*}
\begin{equation}
\begin{aligned}
z_1(x)&=\cosh(\gamma x)\cos(\gamma x),\\
z_2(x)&=\frac{\cosh(\gamma x) \sin(\gamma x)+\sinh(\gamma x)\cos(\gamma x)}{2\gamma},\\
z_3(x)&=\frac{\sinh(\gamma x) \sin(\gamma x)}{2\gamma^2},\\
z_4(x)&=\frac{\cosh(\gamma x) \sin(\gamma x)-\sinh(\gamma x)\cos(\gamma x)}{4 \gamma^3}.
\end{aligned}
\label{z_functions}
\end{equation}

Formula~\eqref{GenSol} represents the solutions of~\eqref{ODE_vector} in terms of their boundary values $\bar W^0$ and $\bar W^l$ at $x=0$ and $x=l$, respectively,
From the boundary conditions~\eqref{BC}, we conclude that
\begin{equation}
\bar W^0=(0,\bar W^0_2,0,\bar W^0_4)^T,\; \bar W^l=(0,\bar W^l_2,0,\bar W^l_4)^T.
\label{boundary_constants}
\end{equation}
To eliminate the parameters $\bar W^0_2$, $\bar W^0_4$, $\bar W^l_2$,
$\bar W^l_4$, we exploit the property that $W^s(x)$ is of class $C^2[0,l]$ together with the interface condition~\eqref{interface}. As a result, we get the following linear algebraic system with respect to $\bar W=(\bar W^0_2, \bar W^0_4, \bar W^l_2,
\bar W^l_4)^T$:
\begin{equation}\label{C_eq}
M \bar W  = R,
\end{equation}
where

\begin{table*}
	\vspace{-\baselineskip}
	\begin{align}
	{\scriptsize
		M= \begin{pmatrix}
		-z_2(l_0) & -z_4(l_0) & z_2(l_0-l) & z_4(l_0-l) \\
		-z_1(l_0) & -z_3(l_0) & z_1(l_0-l) & z_3(l_0-l) \\
		4 \gamma^4 z_4(l_0) & -z_2(l_0) & -4 \gamma^4 z_4(l_0-l) & z_2(l_0-l) \\
		\beta z_2(l_0) + 4 \gamma^4 z_3(l_0)& \beta z_4(l_0) - z_1(l_0) & -4 \gamma^4 z_3(l_0-l) & z_1(l_0-l)
		\end{pmatrix}, R = \begin{pmatrix}
		\int_0^l z_4(l_0-y) \phi^s(y) dy \\
		\int_0^l z_3(l_0-y) \phi^s(y) dy \\
		\int_0^l z_2(l_0-y) \phi^s(y) dy \\
		\frac{U_0}{EI} + \int_0^l z_1(l_0-y) \phi^s(y) dy - \beta \int_0^{l_0} z_4(l_0-y) \phi^s(y)dy
		\end{pmatrix},}\label{MR}
	\end{align}
\end{table*}
and
$$
\beta = \frac{s^2 m+\varkappa}{EI}.
$$
Thus, the vector-valued function $W^s(x)$ is defined by~\eqref{GenSol} and~\eqref{boundary_constants} with $\bar W=M^{-1} R$, and
the components of $R=(R_1^s,R_2^s,R_3^s,R_4^s)^T$ are linear combinations of $U_0$, ..., $U_k$:
\begin{equation}
R_i^s = \sum_{j=0}^k r_{ij}^s U_j,\;\; i=1,2,3,4,
\label{Ri}
\end{equation}
where the coefficient matrix $r^s=(r_{ij}^s)$ can be obtained from~\eqref{phi},~\eqref{MR}.
Then the function $W_3^s(x)$, corresponding to the second $x$-derivative of $w(x,t)$, is expressed as:
$$
\begin{aligned}
& W_3^s(x) = -4\gamma^4 z_4(x)\bar W^0_2 + z_2(x)\bar W_4^0 \\ & + \sum_{j=1}^k \frac{U_j}{EI+ \rho d s } \int_0^x K_j(x,y) dy\;\;\text{for}\; x\in [0,l_0],
\end{aligned}
$$
\begin{equation}
\begin{aligned}
& W_3^s(x) = -4\gamma^4 z_4(x-l)\bar W^l_2 + z_2(x-l)\bar W_4^l \\ & - \sum_{j=1}^k \frac{U_j}{EI+ \rho d s } \int_x^l K_j(x,y) dy\;\text{for}\; x\in (l_0,l],
\end{aligned}\label{W3}
\end{equation}
where $K_j(x,y)=z_2(x-y)\psi_j''(x)$ and $(\bar W^0_2, \bar W^0_4, \bar W^l_2,
\bar W^l_4)^T=M^{-1}r^s (U_0,U_1,...,U_k)^T$ because of~\eqref{Ri}.
Thus, at each $x\in [0,l]$, the above formulas define $W_3^s(x)$ as a linear combination of $U_0$, $U_1$,..., $U_k$:
\begin{equation}
W^s_3(x) = \sum_{j=0}^k h_j^s(x)U_j,
\label{W3h}
\end{equation}
with the coefficients $h_j^s(x)$ collected from~\eqref{W3}.
Recalling that the output of the considered system is given by~\eqref{outputs}, we summarize the computation of the transfer function in the following lemma.
\begin{lemma}\label{lem_TF1}
	The transfer function of the multi-input multi-output control system~\eqref{EulerBernoulli}--\eqref{outputs} is presented in the form:
	$$
	H(s) = \begin{pmatrix}
	h_0^s(l_1) & h_1^s(l_1) & ... & h_k^s(l_1)\\
	h_0^s(l_2) & h_1^s(l_2) & ... & h_k^s(l_2)\\
	\vdots & \vdots & \ddots & \vdots  \\
	h_0^s(l_p) & h_1^s(l_p) & ... & h_k^s(l_p)\\
	\end{pmatrix},
	$$
	where $h_j^s(x)$ are taken from~\eqref{W3h}.
\end{lemma}

\subsection{Single-input single-output (SISO) case}

Let the system be controlled by the shaker force $u_0$ only and the scalar output signal $y_1(t)$ be available.
In this particular case, the scalar transfer function $H(s)$ is such that
$
Y_1(s) = H(s) U_0(s)
$.
Lemma~\ref{lem_TF1} implies the following result in the considered SISO case.
\begin{lemma}\label{lem_TF2}
	Assume that $k=0$, $p=1$, and $l_1\le l_0$. Then the transfer function of the control system~\eqref{EulerBernoulli}--\eqref{outputs} is
	\begin{equation}
	H(s) = \frac{1}{EI}\left( -4\gamma^4 z_4(l_j) M^{-1}_{14} + z_2(l_j) M^{-1}_{24}\right).
	\label{tf_SISO}
	\end{equation}
	Here $M^{-1}_{ik}$ are elements of $M^{-1}$ (the matrix $M$ is given in~\eqref{MR}) and $z_i(x)$ are defined in~\eqref{z_functions}.
\end{lemma}

For further numerical simulations, we
consider the beam actuated by the shaker force only ($k=0$) with single output ($p=1$) and take the following realistic mechanical parameters~\cite{KZB2021} (see also~\cite{D2014}):
$$
l=1.905\,\text{m},\; l_0=1.4\,\text{m},\;\rho_0 = 2700 \,\text{kg}/\text{m}^3,
$$
$$
S = 2.25\cdot10^{-4} \text{m}^2,\; \rho=\rho_0 S,\;
E=6.9\cdot10^{10}\,\text{Pa},
$$
$$
I=1.6875\cdot10^{-10}\,\text{m}^4,\; m=0.1\,\text{kg},\; \varkappa=7\,\text{N/mm},
$$
\begin{equation}
l_1 = 732.5\,\text{mm}.
\label{mechparams}
\end{equation}

\section{The Loewner framework for fitting linear time-invariant systems}\label{sec:loew}

In what follows, we provide a brief summary of the LF for fitting linear dynamical systems as in \eqref{eq:linsys}, from data. The starting point for LF is having access to measurements corresponding to the transfer function of the underlying dynamical process, which can be inferred in practice by means of experimental or model-based procedures. The data set is given by:
\begin{equation}\label{eq:dataset}
\cD = \{(\omega_\ell;H(\omega_\ell))\ \vert
\ \ell=1,\ldots,2k\},
\end{equation}
by means of sampling $H: \IC \rightarrow \IC$ is an analytic function (not necessarily rational) on a particular (complex) grid of points $\omega_\ell$'s. It is to be noted that data sets with an odd number of measurements can also be accommodated in the LF. The first step is to partition the data set in (\ref{eq:dataset}) into two disjoint subsets, as follows:
\begin{align}\label{data_Loew}
\begin{split}
{\textrm right \ data}&: \ \cD_R = \{(\lambda_j;w_j)\ \vert
\ j=1,\ldots,k\},~{\textrm and}, \\
{\textrm left \ data}&: \ \cD_L = \{(\mu_i;v_i) \ \vert \ ~i=1,\ldots,k\},
\end{split}
\end{align}
For simplicity, all points are assumed distinct and also $\mu_i \neq \lambda_j$, for all $1 \leq i, j \leq k$; extensions to Hermite interpolation were proposed in \cite{MA07}.

A typical approach for splitting the data, commonly used in the LF publications, is the ``alternate splitting scheme'', described as follows. The left and right sample nodes (and points) are chosen so that they are interlacing each other. More precisely, for $1 \leq i \leq k$, we can write:
\begin{align}
\begin{split}
\mu_i &= \omega_{2i-1}, \ \ \lambda_i = \omega_{2i},\\
v_i = H(\mu_i) = H(&\omega_{2i-1}), \ \ w_i = H(\lambda_i) = H(\omega_{2i}).
\end{split}
\end{align}
We refer the reader to \cite{karachalios2021loewner} for a more comprehensive account of data partitioning strategies in the Loewner framework; there, the ``half-hal'' splitting is mentioned together with approaches that split the data based on the magnitude of the data samples.

The goal is to find a rational function denoted with
$\tilde{H}(s)$, such that the following interpolation conditions are (approximately) fulfilled:
\begin{equation} \label{interp_cond}
\tilde{H}(\mu_i)=v_i,~~~\tilde{H}(\lambda_j)=w_j.
\end{equation}
In order to accomplish this scope, we first arrange the elements of the original data set $\cD$, partitioned as in (\ref{data_Loew}) in matrix format. Hence, the Loewner matrix $\IL \in\IC^{k\times k}$ and the shifted Loewner matrix $\IM \in\IC^{k\times k}$ are defined as follows
\begin{equation} \label{Loew_mat}
\IL_{(i,j)}=\frac{v_i-w_j}{\mu_i-\lambda_j}, \ \IM_{(i,j)}=
\frac{\mu_i v_i-\lambda_j w_j}{\mu_i-\lambda_j},
\end{equation}
while the data vectors $\IV, \IW^T \in \IR^k$ are given by:
\begin{equation} \label{VW_vec}
\IV_{(i)}= v_i, \ \  \IW_{(j)} = w_j,~\text{for}~i,j=1,\ldots,k.
\end{equation}
The Loewner model is hence constructed as follows:
\begin{align*}
\bE=-\IL,~~ \bA=-\IM,~~ \bB=\IV,~~ \bC=\IW.
\end{align*}
The following Sylvester equations are satisfied by the Loewner and shifted Loewner matrices, as shown in \cite{ALI17} (here, $\bone_q = \left[ \begin{matrix} 1 &  \cdots & 1  \end{matrix} \right]^T \in \IC^q$):
\vspace{-2mm}
\begin{equation}\label{eq:sylv_loe}
\begin{cases} \bM \IL - \IL \bLambda = \IV \bone_k^T - \bone_q \IW, \\
\bM \IM - \IM \bLambda =  \bM \IV \bone_k^T -  \bone_q  \IW \bLambda,
\end{cases}
\vspace{-2mm}
\end{equation}
where $\bM = \text{diag}(\mu_1,\cdots,\mu_q)$ and $\bLambda = \text{diag}(\lambda_1,\cdots,\lambda_k)$. The following relations expressing the shifted Loewner matrix to the Loewner matrix, in two distinct ways, hold:
\begin{equation}
\IM =  \IL \bLambda + \IV \bone_k^T = \bM \IL + \bone_q  \IW.
\end{equation}
Hence, the explicit computation of large Loewner matrices can be avoided by means of computing (approximated, low-rank) solutions of the Sylvester equations in (\ref{eq:sylv_loe}). This can be accomplished, e.g., by means of using optimized and robust numerical tools such as \cite{morBenKS21}.

Provided that enough data are available, the pencil $(\IM,\,\IL)$ is often singular. For example, if the data $\cD$ in \ref{eq:dataset} were generated from a rational function $H(s)$ corresponding to a minimal LTI system of dimension $n$, i.e., $H(s) = \bC (s\bI_n-\bA)^{-1} \bB$, then this corresponds to the case $k > n$. Then, in the perfect setup (no noisy data), one would encounter $k-n$ zero singular values when performing the singular value decomposition (SVD) of the pencil $\zeta \IL - \IM$, where $\zeta \in \IC$ is chosen to be different than the eigenvalues of matrix $\bA$. In such cases, an SVD of augmented Loewner matrices is computed, and the dominating part is selected as:
\begin{equation}\label{eq:AugmLoew}
\left[\IL,~\IM\right] \approx \bY\widehat{\Sigma}_{ {r}}\tilde{\bX}^*,~\left[\begin{array}{c}\IL \\ \IM\end{array}\right] \approx{\tilde\bY}\Sigma_{ {r}} \bX^*,
\end{equation}
where $\widehat{\Sigma}_{ {r}}$, $\Sigma_{ {r}}$ $\in\IR^{{{r}}\times{r}}$,~
$\bY \in\IC^{k\times{r}}$,$\bX\in\IC^{k\times{r}}$,~$\tilde{\bY}\in\IC^{2k\times{r}}$,~$\tilde{\bX}\in\IC^{r\times{2k}}$ and $(\bX)^* \in \IC^{r \times k}$ denotes the conjugate-transpose of matrix $\bX$. This is performed in order to find projection matrices $\bX_r, \bY_r \in \IC^{k \times r}$, as described in \cite{ALI17}. Here, $r<n$ represents the truncation index. Then, the system matrices corresponding to a projected Loewner model of dimension $r$  can be computed using the truncated singular vector matrices $\bX_r$ and $\bY_r$:
\begin{align}\label{Loew_red_lin}
\begin{split}
\hat{\bE} &= -\bX_r^*\IL \bY_r, \ \  \hat{\bA} = -\bX_r^*\IM \bY_r, \\
\hat{\bB} &= \bX_r^*\IV, \ \  \hat{\bC} = \IW \bY_r,
\end{split}
\end{align}
and therefore, directly finds a state-space realization corresponding to the reduced-order system of equations
\begin{equation}\label{eq:linsys_red}
\begin{cases}
\bhE \dot{\bhx}(t)=\bhA\bx(t)+\bhB u(t),\\ \hat{y}(t) \ \hspace{1.8mm}=\bhC\bhx(t).
\end{cases}.
\end{equation}
The transfer function of the reduced Loewner model in (\ref{eq:linsys_red}) is written as $\hH(s) = \bhC (s \bhE - \bhA)^{-1} \bhB$, and it provides a good approximant to the original transfer function $H(s)$. Then, $\hH(s)$ may be expanded in a pole/zero or pole/residue format. These values represent system invariants and can be related to the inherent dynamics. It is noted that the state-space realization is not unique, and that is why an extra step is required. More implementation details and properties of the LF procedure can be found in \cite{ALI17,karachalios2021loewner}.

\section{Numerical examples}

\subsection{Analysis on the unperturbed data}

In this section, we present two numerical test cases based on sampling the transfer function explicitly derived in Section~\ref{sec_tf}.
We consider the SISO system with the choice of physical parameters~\eqref{mechparams}.
Then the scalar transfer function $H(s)$ in Lemma~\ref{lem_TF2} is sampled at the purely imaginary grid points $s=\omega_\ell$, $\ell=1,\ldots,k$ with $k=1000$, which are equally distributed in the range of physical frequencies from $0 Hz$ to $250 Hz$. The data partitioning scheme (into the left and right disjoint subsets) chosen here is the ``alternate'' one, previously described in Section \ref{sec:loew}. We also note that complex conjugated data is added to the process to enforce real-valued models; more details on how this is achieved can be found in \cite{ALI17}.

\begin{center}
	\begin{figure}[h!]
		\hspace{-6mm}\includegraphics[scale=0.24]{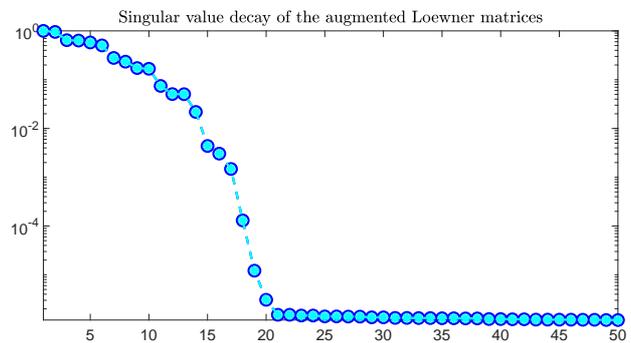}
		\caption{The decay of the singular values for the augmented Loewner matrices in \ref{eq:AugmLoew}.}
		\label{fig1}
	\end{figure}
\end{center}

\begin{figure}[h!]\hspace{-3mm}\includegraphics[scale=0.24]{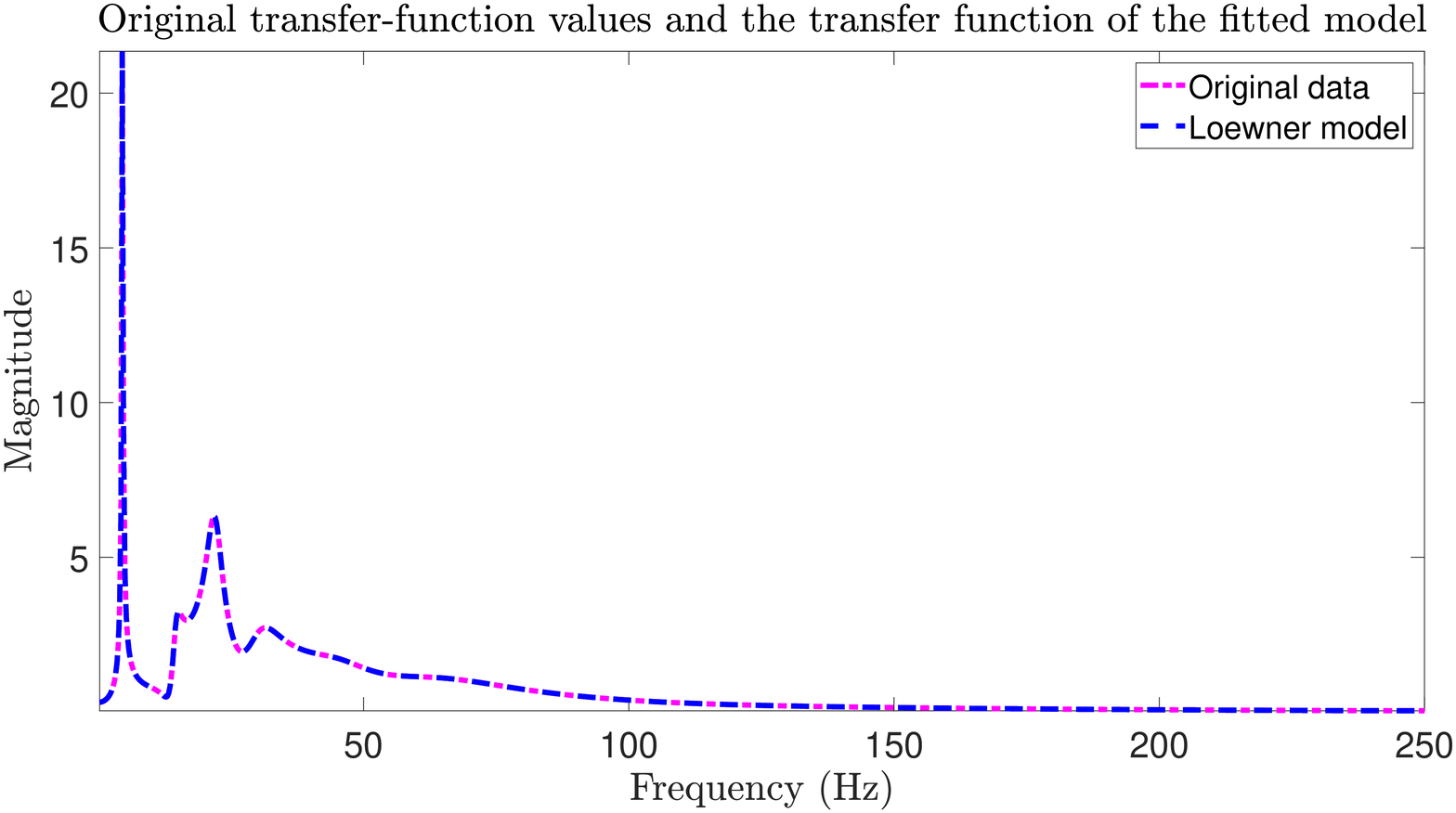}
	\caption{The original data (transfer function samples) vs. the Loewner model fit.}
	\label{fig2}
\end{figure}

\begin{figure}[h!]
	\hspace{-4mm}\includegraphics[scale=0.24]{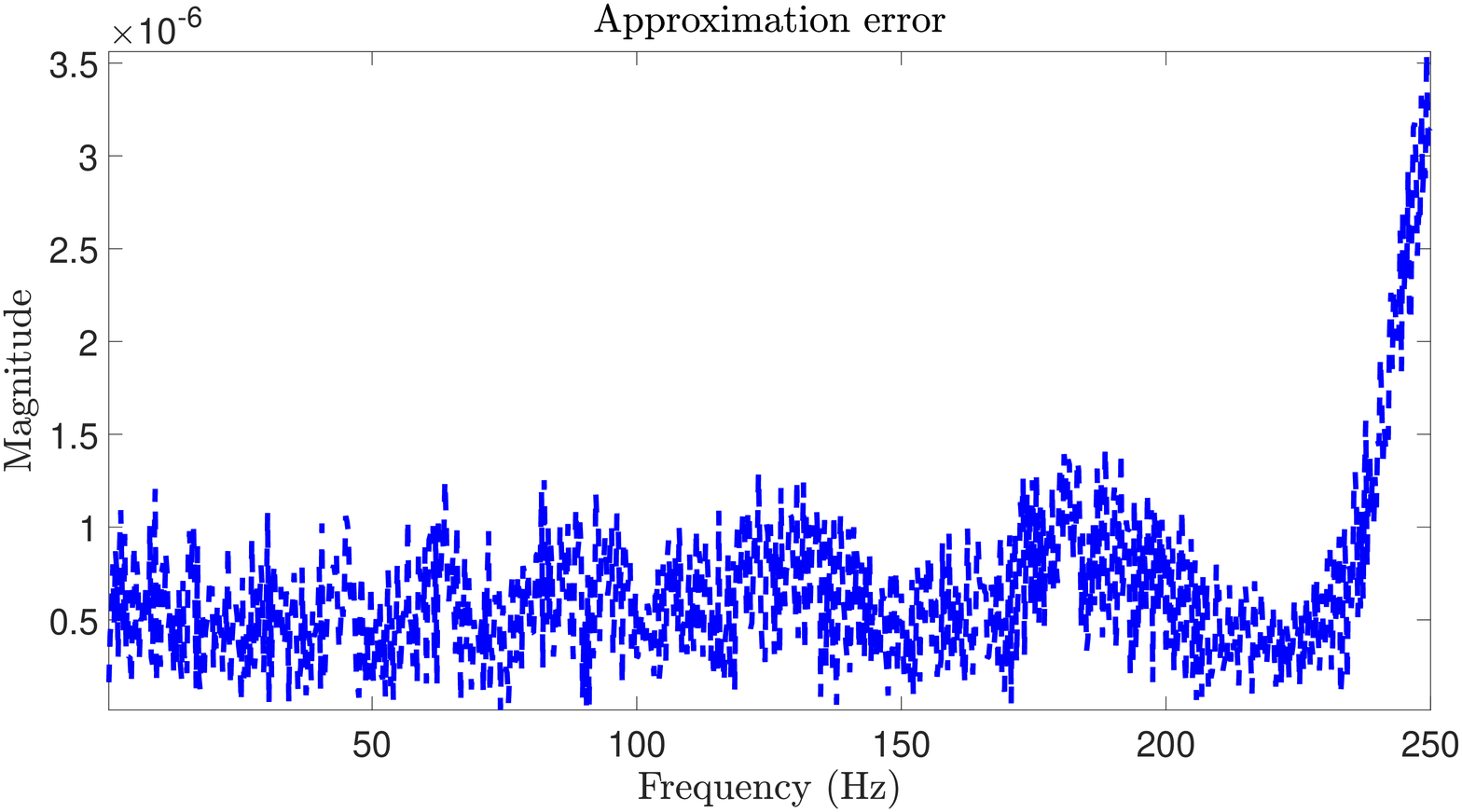}
	\caption{The approximation error.}
	\label{fig3}
\end{figure}

\begin{figure}[h!]
	\hspace{-4mm}\includegraphics[scale=0.24]{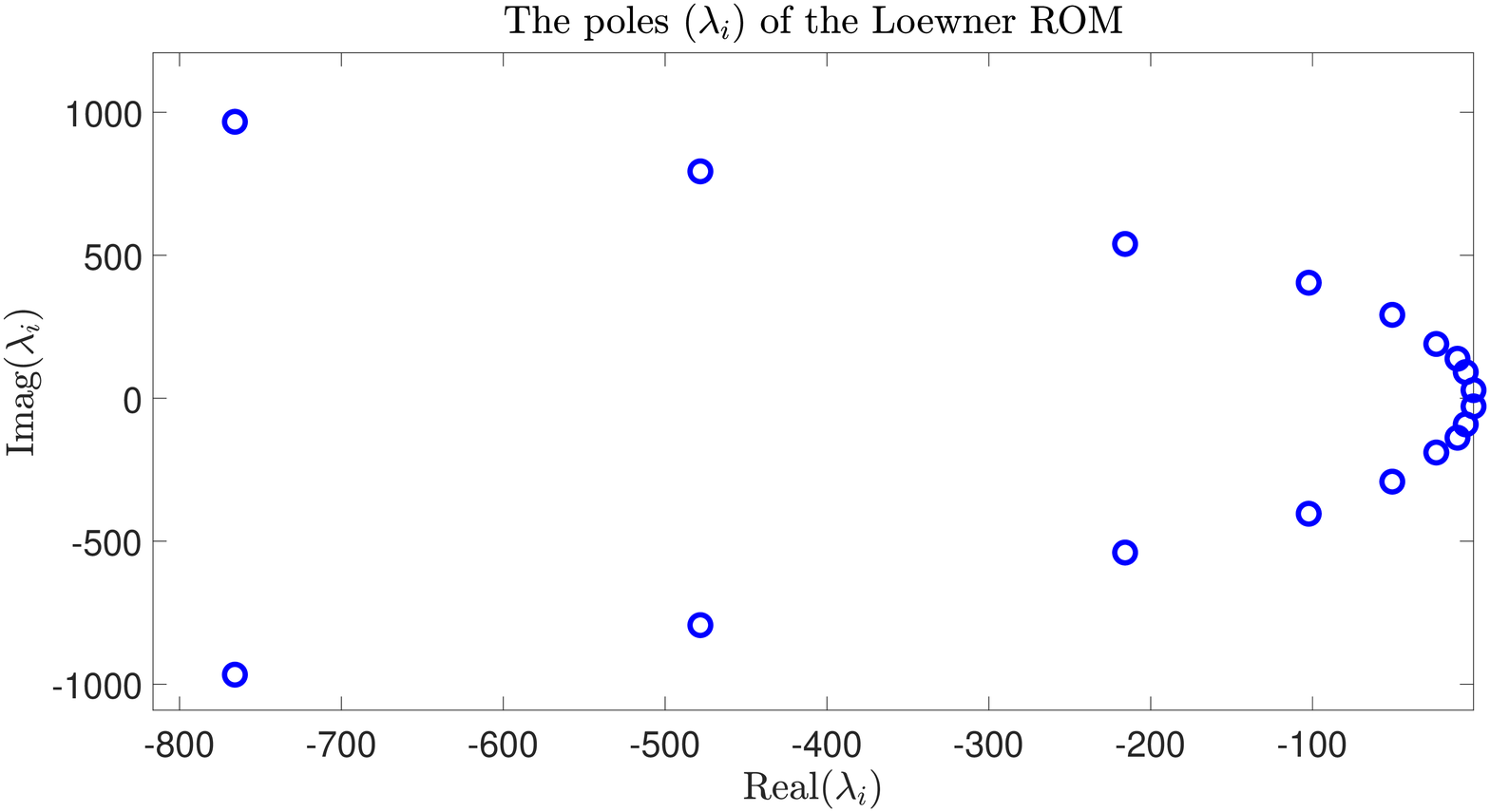}
	\caption{The poles of the fitted Loewner model, as eigenvalues of pencil $(\bhA,\bhE)$.}
	\label{fig4}
\end{figure}

\begin{center}
	\begin{figure}[h!]
		\hspace{-4mm}\includegraphics[scale=0.24]{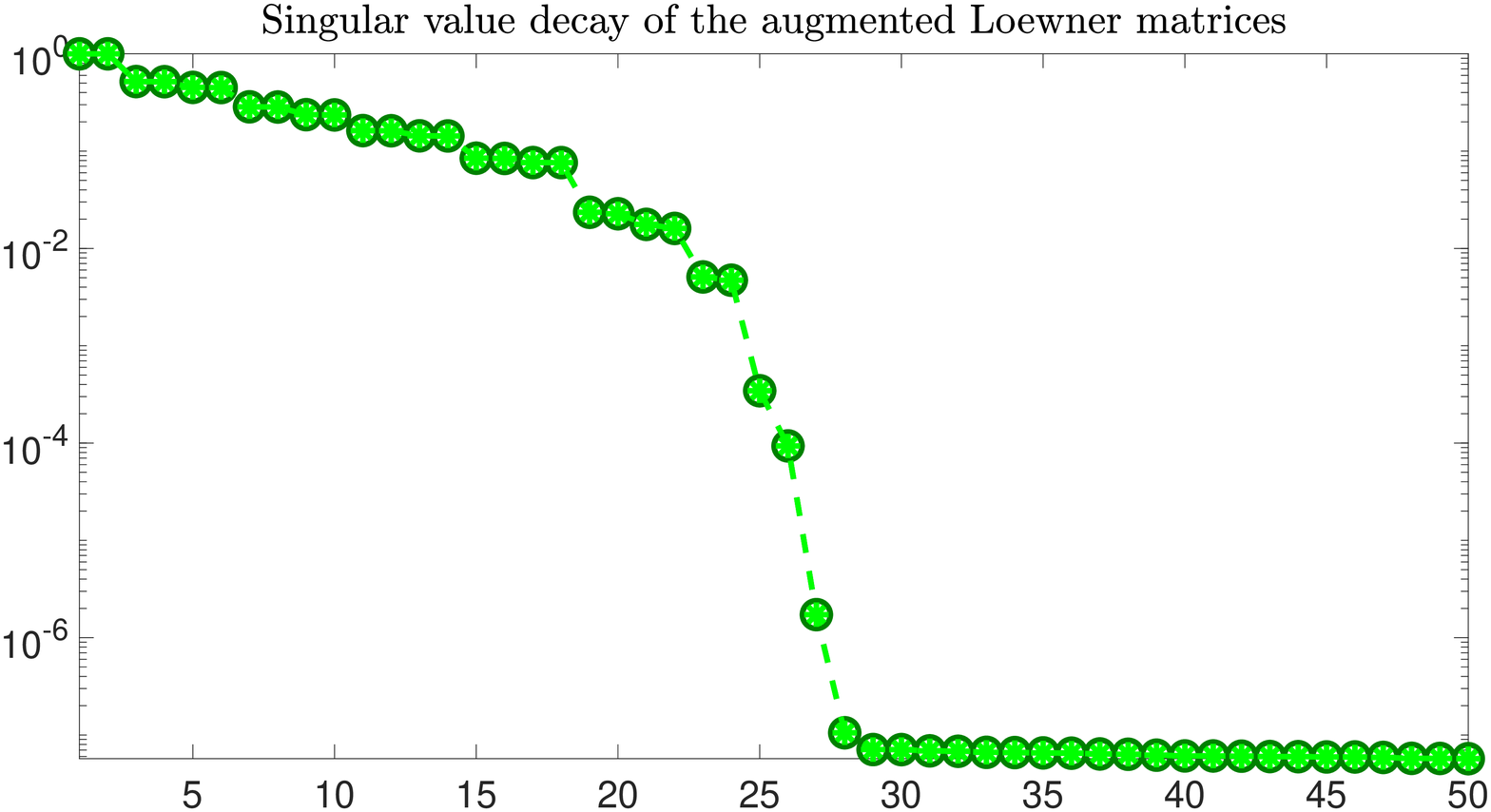}
		\caption{The decay of the singular values for the augmented Loewner matrices in \ref{eq:AugmLoew}.}
		\label{fig5}
	\end{figure}
\end{center}

\begin{figure}[h!]\hspace{-3mm}\includegraphics[scale=0.24]{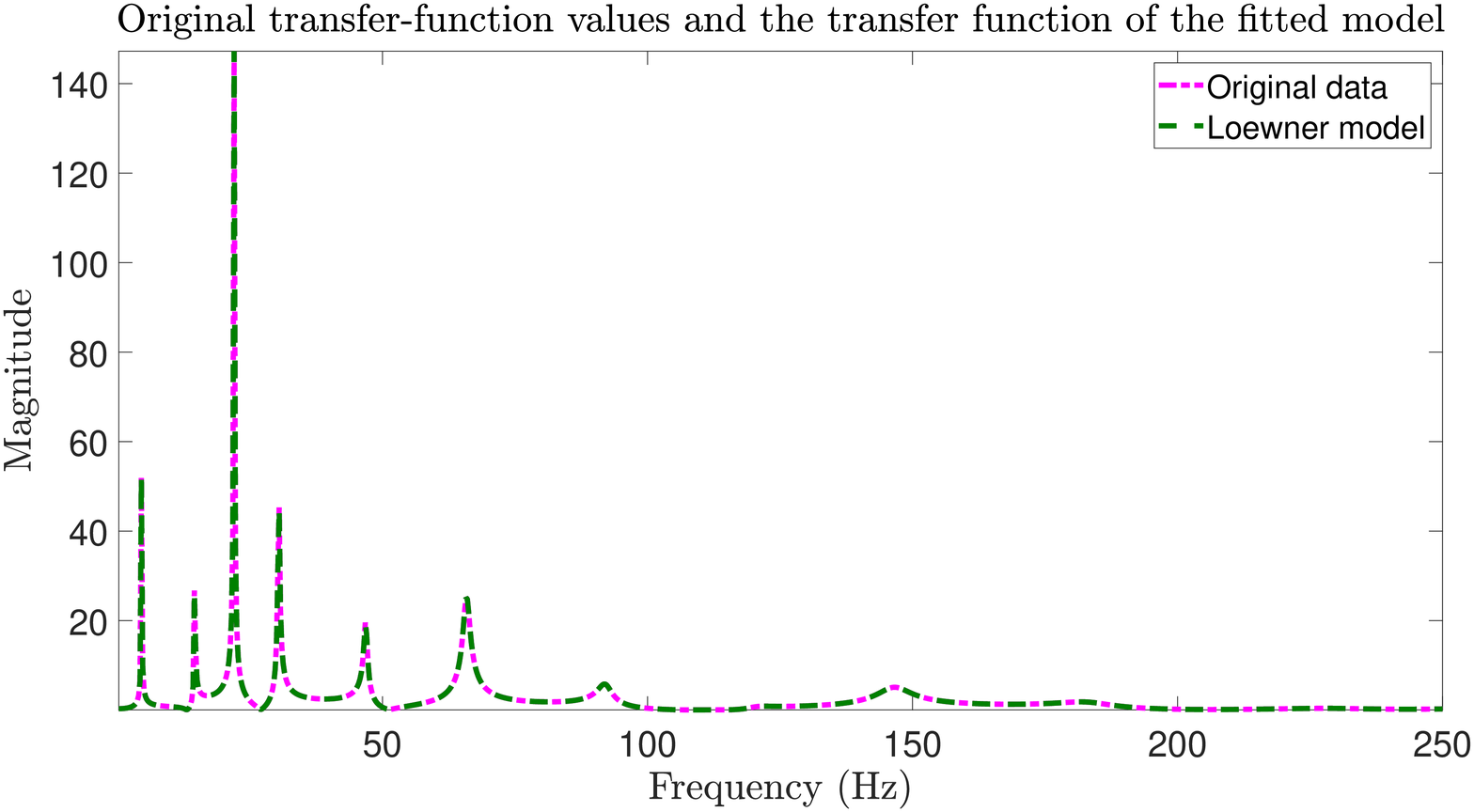}
	\caption{The original data (transfer function samples) vs. the Loewner model fit.}
	\label{fig6}
\end{figure}

\begin{figure}[h!]
	\hspace{-4mm}\includegraphics[scale=0.24]{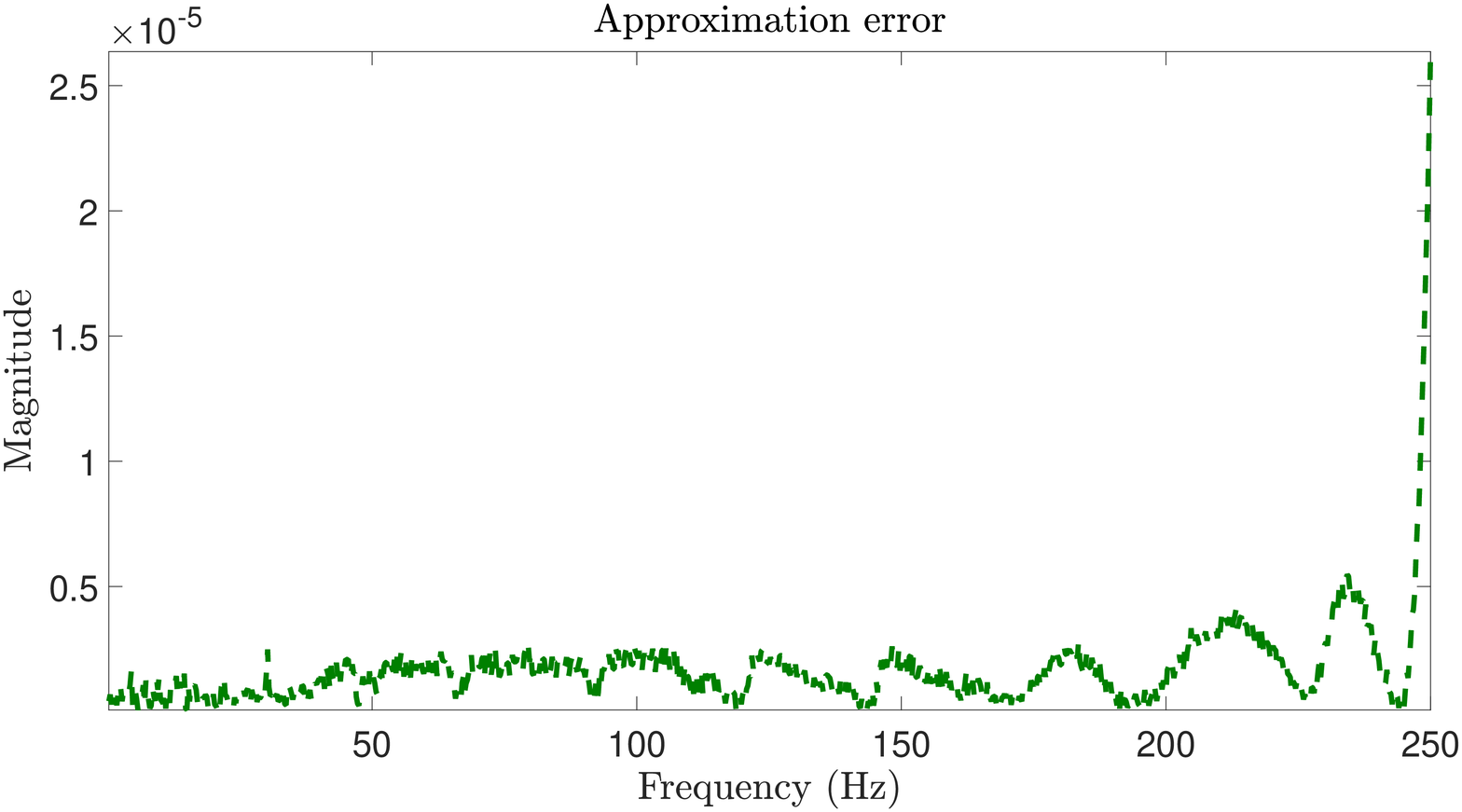}
	\caption{The approximation error.}
	\label{fig7}
\end{figure}

\begin{figure}[h!]
	\hspace{-4mm}\includegraphics[scale=0.24]{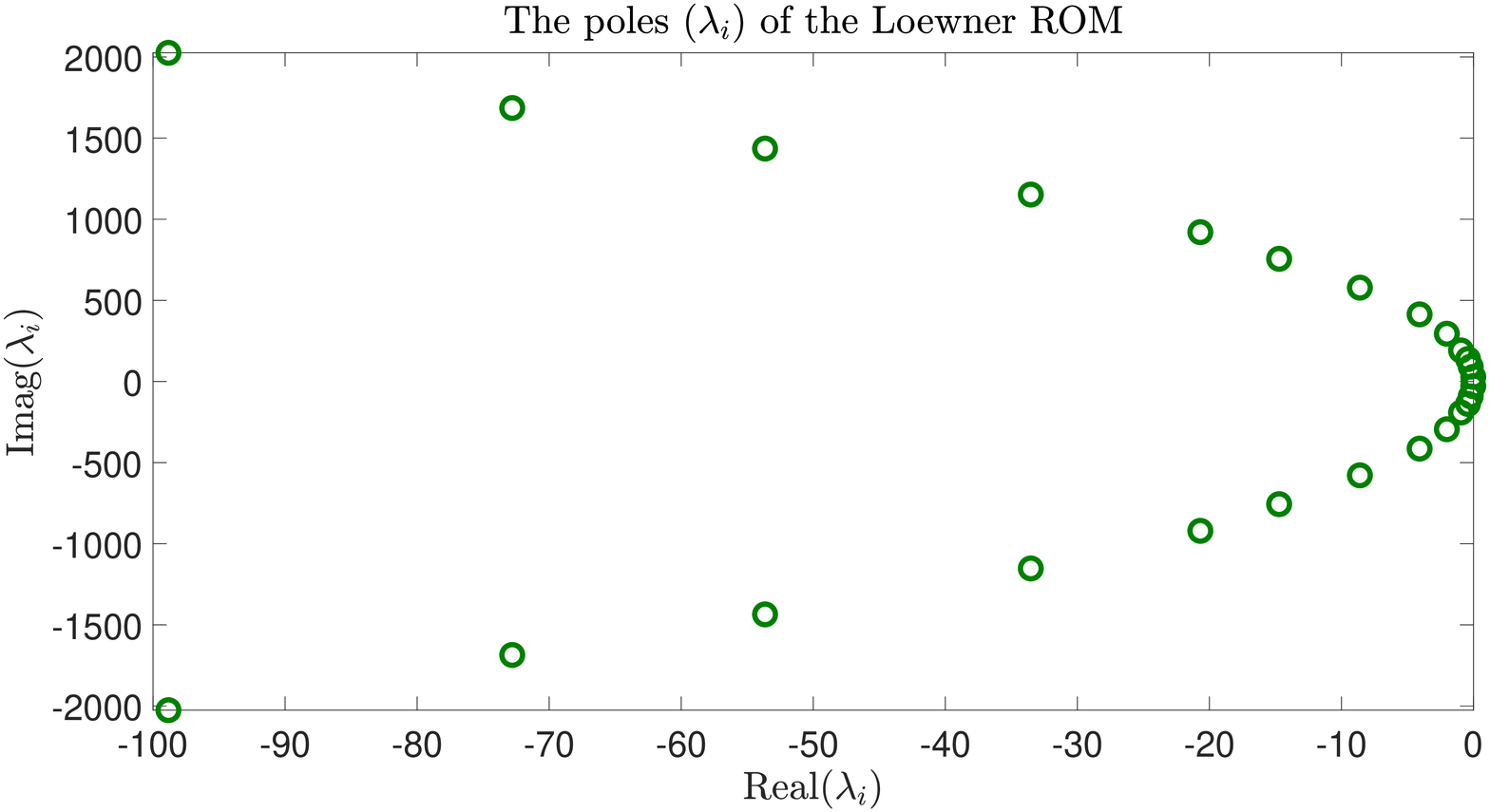}
	\caption{The poles of the fitted Loewner model, as eigenvalues of pencil $(\bhA,\bhE)$.}
	\label{fig8}
\end{figure}

Two values of the structural damping parameter $d$ in~\eqref{EulerBernoulli} are considered for the numerical simulations: $d_1=0.0249$ (``large'' damping) and $d_2=0.001$ (``small'' damping).
The case $d=d_1$ is depicted in Figs.~\ref{fig1}--\ref{fig4}, for which we are fitting a  Loewner model of order $r = 20$ (with a rational transfer function). Additionally, the case $d=d_2$ is presented in Figs.~\ref{fig5}--\ref{fig8}; for this case, we are fitting a  Loewner model of order $r = 27$. It should be noted that the Loewner framework does not automatically impose stability; post-processing methods can be applied whenever unstable poles appear, as described in~\cite{GA16}.

We observe that the Loewner model approximates the original transfer function with good accuracy in the whole range of test frequencies (the approximation error of order $10^{-6}$ in Fig.~\ref{fig3}, and of order $10^{-5}$ in Fig.~\ref{fig7}). It is also clear that these approximate models preserve the stability property (the maximal real part of $\lambda_i$ in Fig.~\ref{fig4} is ${\textrm max}({\textrm Re}\,\lambda_i) = -5.2780  \cdot 10^{-1}<0$, while the maximal real part of $\lambda_i$ in Fig.~\ref{fig8} is ${\textrm max}({\textrm Re}\,\lambda_i) = -2.1187 \cdot 10^{-2}<0$).

Although the transfer function of the considered distributed parameter function is not rational, the decay of singular values associated with the Loewner matrices in Figs.~\ref{fig1} and \ref{fig5} seems to indicate the opportunity to enforce rational approximation. More precisely, in both cases, a plateau (flat portion of the graph) is observed after a steep decay. In Figs.~\ref{fig1} and \ref{fig5}, we only depict the first $50$ singular values (out of $1000$, which is the dimension of the Loewner matrices $\IL$ and $\IM$). In the first test case, the decay is faster than in the second one. The dimensions of the reduced-order models were chosen in accordance with this phenomenon (i.e. $r=20$ for the first and $r=27$ for the second). However, such a clear and steep decay as seen in Figs.~\ref{fig1} and \ref{fig5}, is seldom noticed in experimental data. In such scenarios, the data can be perturbed, i.e., by means of noise. We treat this case below.

\subsection{Analysis on the perturbed data (by means of artificial additive Gaussian noise)}

In this subsection, we analyze the robustness of the LF when applied to perturbed (noisy) data. A preliminary analysis of such endeavors was reported in~\cite{lefteriu2010modeling,drmavc2022learning,kergus2022data} and in \cite{zhang2021factorization, embree2022pseudospectra}. Similarly to the approaches in these publications, we will include additive Gaussian noise into the measurements of the transfer function $H(s)$ from Lemma~\ref{lem_TF2}. More precisely, for all $k=1000$ previous grid points $s=\omega_\ell$, the new data are, for $\ell=1,\ldots,k$ and $\nu = 1,\ldots,4$:
\begin{equation}
H( \omega_\ell)[1+\epsilon^{(\nu)}(\alpha_\ell+\imath \beta_\ell)].
\end{equation}
Here, the ``noise power'' $\epsilon^{(\nu)} >0$ is chosen so that $\epsilon^{(\nu)} = 10^{-\nu}$, while the values $\alpha_\ell$ and $\beta_\ell$ are drawn from the standard normal distribution. For simplicity, we use the term ``noise level $\nu$'' (in the subsequent text and plots). Also, for the sake of brevity, we will restrict our analysis to the case of a smaller damping coefficient: $d = 0.001$.

\begin{figure}[h!]
	\hspace{-4mm}\includegraphics[scale=0.24]{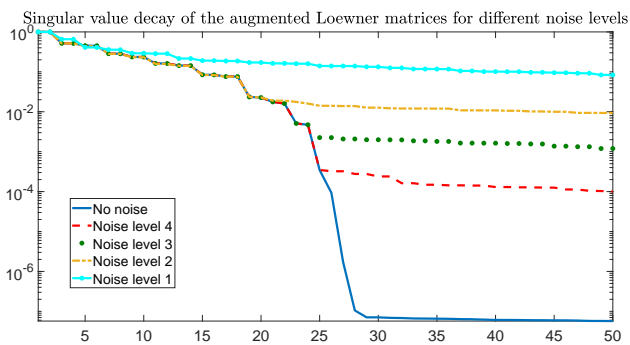}
	\caption{The decay of singular values for the augmented Loewner matrices  for different "noise levels/powers''.}
	\label{fig9}
\end{figure}

The first numerical experiment that is shown here is concerned with the decay of singular values for Loewner matrices under the influence of noise. As previously pointed out in \cite{lefteriu2010modeling,gosea2021data}, the effect of noise is generally reflected in the ``flattening'' of the singular value curve. This is especially valid for the ``alternate'' splitting, and not as much for the ``half-half'' splitting (as shown in \cite{gosea2021data}). This is precisely the phenomenon observed in Fig.~\ref{fig9}. Additionally, there are a number of dominant singular values that are less prone to perturbations; depending on the noise levels, it is clear that this number decreases with the increase of the noise power. For example, when $\nu =2$ there are 22 singular values that seem to be stagnant. This is hence another effect of the noise in data for LF; deciding the order of the fitted model becomes a more challenging task, and one needs to be careful to avoid over-fitting.

\begin{figure}[h!]
	\hspace{-4mm}\includegraphics[scale=0.24]{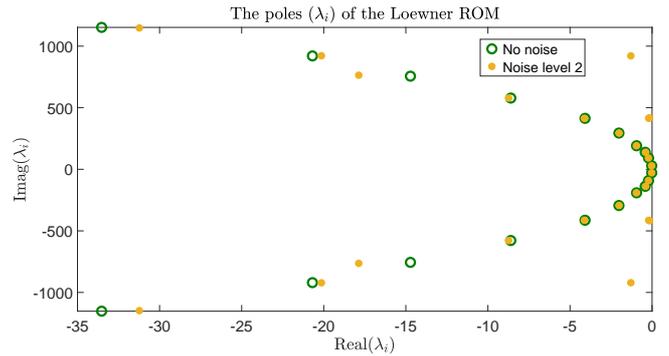}
	\caption{The poles of the fitted Loewner models (with and without noise).}
	\label{fig10}
\end{figure}

Next we fix the noise level at $\nu =2$ and  record the poles of the fitted Loewner model on the noisy data. We compare those with the original poles, and the results are depicted in Fig.~\ref{fig10}. As expected, most of the (dominant) poles seem to be matching well, with the remark that the noisy data introduces spurious poles.

\begin{figure}[h!]
	\hspace{-4mm}\includegraphics[scale=0.24]{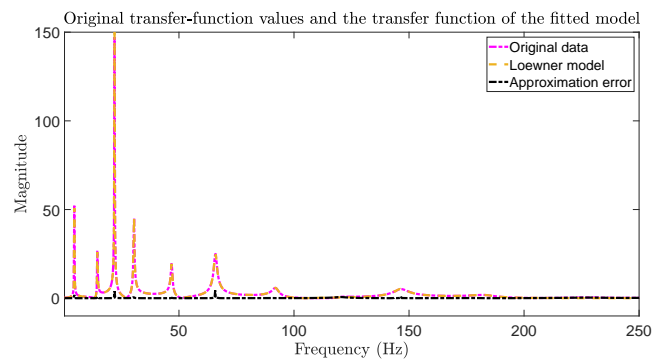}
	\caption{The original data vs. the Loewner model fit (on the noisy data for level $\nu=2$) $\&$ the approximation error.}
	\label{fig11}
\end{figure}

However, as shown in Fig.~\ref{fig11}, the effects of the noise (for this level) do not seem to be drastic; the response of the Loewner model fitted to the noisy data faithfully follows the original (unperturbed) data.

\begin{figure}[h!]
	\hspace{-4mm}\includegraphics[scale=0.24]{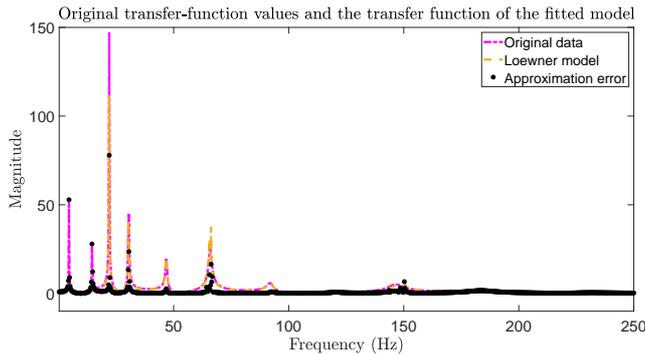}
	\caption{The original data vs. the Loewner model fit (on the noisy data for level $\nu=1$) $\&$ the approximation error.}
	\label{fig12}
\end{figure}

It is to be noted that for levels of noise higher than the one used in the previous experiment, i.e., for $\nu=1$, the results are significantly less accurate. In that case, for some perturbed data sets, the first two dominant peaks are partially or completely missed (as shown in Fig.~\ref{fig12}). Additionally, it was noticed that the asymptotic stability of the reduced-order Loewner model was sometimes lost (however not in the experiment reported here). This behavior is mostly due to the high level of the noise signal added here. As reported in \cite{lefteriu2010modeling}, for high signal-to-noise scenarios, it is very challenging to extract all the meaningful information from the perturbed data. A more thorough numerical analysis based on the so-called stabilization diagrams will be left for future research endeavors.

\section{Conclusion and outlook}


The contribution of this work is twofold. First, an analytic construction of the transfer function for an infinite-dimensional flexible structure with the Euler--Bernoulli beam and a spring-mass system has been proposed.
Second, we have applied the Loewner framework (LF) for data-driven modeling of the considered class of flexible structures on the basis of transfer function measurements. The presented case study of the singular value decay for the Loewner pencil clearly indicates a possible choice of the dimension of an acceptable reduced-order model, depending on the damping parameter $d$. Note that our results are not limited to deterministic measurements since the performed study evaluates the effects of noise in the LF for data sets with Gaussian perturbations as well.
A thorough numerical analysis of the perturbed (noisy) data has been carried out. The preliminary results show the robustness of the method to low and moderate levels of noise, but also point out some challenges for the cases in which the data are perturbed with higher noise levels. For future research endeavors, we intend to tackle the following open issues:
\begin{itemize}
	\item Analyze data from real experimental measurements to supplement the information attained from simulation data.
	\item Take into account and quantify the effects of measurement noise in the LF for such data sets (using pseudospectra theory \cite{embree2022pseudospectra}).
	\item Study the decay of the singular values depending on the influence of the damping parameter and on data splitting strategies. Extend this study to flexible structures with a different asymptotic distribution of the eigenvalues, e.g., the Timoshenko beam with attached rigid bodies~\cite{zuyev2007stabilization}.
	\item Examine the applicability of the Loewner-based reduced-order systems for the control design of the original infinite-dimensional plant.
\end{itemize}



\bibliographystyle{plain}
\bibliography{IVGref}

\begin{thebibliography}{10}

\bibitem{altiner2019modeling}
B.~Alt{\i}ner, A.~Deliba{\c{s}}{\i}, and B.~Erol.
\newblock Modeling and control of flexible link manipulators for unmodeled
  dynamics effect.
\newblock {\em Proceedings of the Institution of Mechanical Engineers, Part I:
  Journal of Systems and Control Engineering}, 233(3):245--263, 2019.

\bibitem{ACA05}
A.~C. Antoulas.
\newblock {\em Approximation of large-scale dynamical systems}.
\newblock SIAM, Philadelphia, 2005.

\bibitem{ABG20}
A.~C. Antoulas, C.~A. Beattie, and S.~Gugercin.
\newblock {\em Interpolatory Methods for Model Reduction}.
\newblock SIAM, Philadelphia, 2020.

\bibitem{ALI17}
A.~C. Antoulas, S.~Lefteriu, and A.~C. Ionita.
\newblock A tutorial introduction to the {L}oewner framework for model
  reduction.
\newblock In {\em Model Reduction and Approximation}, chapter~8, pages
  335--376. SIAM, 2017.

\bibitem{beattie2009interpolatory}
C.~Beattie and S.~Gugercin.
\newblock Interpolatory projection methods for structure-preserving model
  reduction.
\newblock {\em Systems \& Control Letters}, 58(3):225--232, 2009.

\bibitem{morBenGKetal20}
P.~Benner, P.~Goyal, B.~Kramer, B.~Peherstorfer, and K.~Willcox.
\newblock Operator inference for non-intrusive model reduction of systems with
  non-polynomial nonlinear terms.
\newblock {\em Comp. Methods in App. Mechanics and Engineering}, 372, 2020.

\bibitem{benner2021structure}
P.~Benner, S.~Gugercin, and S.~W.~R. Werner.
\newblock Structure-preserving interpolation for model reduction of parametric
  bilinear systems.
\newblock {\em Automatica}, 132:109799, 2021.

\bibitem{morBenKS21}
P.~Benner, M.~K{\"o}hler, and J.~Saak.
\newblock Matrix equations, sparse solvers: {M-M.E.S.S.}-2.0.1 -- philosophy,
  features and application for (parametric) model order reduction.
\newblock In {\em Model Reduction of Complex Dynamical Systems}, volume 171 of
  {\em International Series of Numerical Mathematics}, pages 369--392.
  Birkh{\"a}user, Cham, 2021.

\bibitem{BOCW17}
P.~Benner, M.~Ohlberger, A.~Cohen, and K.~Willcox.
\newblock {\em Model Reduction and Approximation}.
\newblock Society for Industrial and Applied Mathematics, Philadelphia, PA,
  2017.

\bibitem{casagrande2019integer}
D.~Casagrande, W.~Krajewski, and U.~Viaro.
\newblock The integer--order approximation of fractional--order systems in the
  Loewner framework.
\newblock {\em IFAC-PapersOnLine}, 52(3):43--48, 2019.

\bibitem{curtain2009transfer}
R.~Curtain and K.~Morris.
\newblock Transfer functions of distributed parameter systems: A tutorial.
\newblock {\em Automatica}, 45(5):1101--1116, 2009.

\bibitem{drmavc2022learning}
Z.~Drma{\v{c}} and B.~Peherstorfer.
\newblock Learning low-dimensional dynamical-system models from noisy
  frequency-response data with Loewner rational interpolation.
\newblock In {\em Realization and Model Reduction of Dynamical Systems}, pages
  39--57. Springer, 2022.

\bibitem{D2014}
C.~Dullinger, A.~Schirrer, and M.~Kozek.
\newblock Advanced control education: optimal \& robust MIMO control of a
  flexible beam setup.
\newblock {\em IFAC Proceedings Volumes}, 47(3):9019--9025, 2014.

\bibitem{embree2022pseudospectra}
M.~Embree and A.~C. Ionita.
\newblock Pseudospectra of {L}oewner matrix pencils.
\newblock In {\em Realization and Model Reduction of Dynamical Systems}, pages
  59--78. Springer, 2022.

\bibitem{GA16}
I.~V. Gosea and A.~C. Antoulas.
\newblock Stability preserving post-processing methods applied in the Loewner
  framework.
\newblock In {\em 2016 IEEE 20th Workshop on Signal and Power Integrity (SPI)},
  pages 1--4, 2016.

\bibitem{gosea2022iterative}
I.~V. Gosea and I.~Pontes~Duff.
\newblock An iterative realization-free approach for model reduction of
  bilinear systems via {H}ermitian interpolation.
\newblock In {\em 2022 European Control Conference (ECC)}, pages 584--589.
  IEEE, 2022.

\bibitem{gosea2021loewner}
I.~V. Gosea, C.~Poussot-Vassal, and A.~C. Antoulas.
\newblock On {L}oewner data-driven control for infinite-dimensional systems.
\newblock In {\em 2021 European Control Conference (ECC)}, pages 93--99. IEEE,
  2021.

\bibitem{gosea2021data}
I.~V. Gosea, Q.~Zhang, and A.~C. Antoulas.
\newblock Data-driven modeling from noisy measurements.
\newblock {\em PAMM}, 20(S1):e202000358, 2021.

\bibitem{morGusS99}
B.~Gustavsen and A.~Semlyen.
\newblock Rational approximation of frequency domain responses by vector
  fitting.
\newblock {\em IEEETransPD}, 14(3):1052--1061, 1999.

\bibitem{KZ2021}
J.~Kalosha and A.~Zuyev.
\newblock Asymptotic stabilization of a flexible beam with an attached mass.
\newblock {\em Ukrainian Mathematical Journal}, 73:1537--1550, 2022.

\bibitem{KZB2021}
J.~Kalosha, A.~Zuyev, and P.~Benner.
\newblock On the eigenvalue distribution for a beam with attached masses.
\newblock In {\em Stabilization of Distributed Parameter Systems: Design
  Methods and Applications}, pages 43--56. Springer, 2021.

\bibitem{karachalios2021loewner}
D.~S. Karachalios, I.~V. Gosea, and A.~C. Antoulas.
\newblock The {L}oewner framework for system identification and reduction.
\newblock In {\em Model Order Reduction: Volume I: System-and Data-Driven
  Methods and Algorithms}, pages 181--228. De Gruyter, 2021.

\bibitem{kergus2022data}
P.~Kergus and I.~V. Gosea.
\newblock Data-driven approximation and reduction from noisy data in matrix
  pencil frameworks.
\newblock {\em arXiv preprint arXiv:2202.09568, to appear in the MTNS22
  proceedings volume}, 2022.

\bibitem{kutz2016dynamic}
J.~N. Kutz, S.~L. Brunton, B.~W. Brunton, and J.~L. Proctor.
\newblock {\em Dynamic mode decomposition: data-driven modeling of complex
  systems}.
\newblock SIAM, 2016.

\bibitem{lefteriu2010modeling}
S.~Lefteriu, A.~C. Ionita, and A.~C. Antoulas.
\newblock Modeling systems based on noisy frequency and time domain
  measurements.
\newblock In {\em Perspectives in Mathematical System Theory, Control, and
  Signal Processing}, pages 365--378. Springer, 2010.

\bibitem{LM18}
P.~Lietaert and K.~Meerbergen.
\newblock Comparing {L}oewner and {K}rylov based model order reduction for time
  delay systems.
\newblock In {\em 2018 European Control Conference (ECC)}, pages 545--550,
  2018.

\bibitem{MA07}
A.~J. Mayo and A.~C. Antoulas.
\newblock A framework for the solution of the generalized realization problem.
\newblock {\em Linear Algebra and Its Applications}, 425(2-3):634--662, 2007.

\bibitem{NST18}
Y.~Nakatsukasa, O.~Sete, and L.~N. Trefethen.
\newblock The {AAA} algorithm for rational approximation.
\newblock {\em SIAM Journal on Scientific Computing}, 40(3):A1494--A1522, 2018.

\bibitem{morPehW16}
B.~Peherstorfer and K.~Willcox.
\newblock Data-driven operator inference for nonintrusive projection-based
  model reduction.
\newblock {\em Computer Methods in Applied Mechanics and Engineering},
  306:196--215, 2016.

\bibitem{duff2015realization}
I.~Pontes-Duff, C.~Poussot-Vassal, and C.~Seren.
\newblock Realization independent single time-delay dynamical model
  interpolation and $\mathcal{H}_2$-optimal approximation.
\newblock In {\em 2015 IEEE Conference on Decision and Control (CDC)}, pages
  4662--4667. IEEE, 2015.

\bibitem{poussot2022interpolation}
C.~Poussot-Vassal, P.~Kergus, and P.~Vuillemin.
\newblock Interpolation-based irrational model control design and stability
  analysis.
\newblock In {\em Realization and Model Reduction of Dynamical Systems}, pages
  353--371. Springer, 2022.

\bibitem{morQuaR14}
A.~Quarteroni and G.~Rozza.
\newblock {\em Reduced Order Methods for Modeling and Computational Reduction},
  volume~9 of {\em MS{\&}A -- Modeling, Simulation and Applications}.
\newblock Springer International Publishing, Cham, Switzerland, 2014.

\bibitem{schulze2016data}
P.~Schulze and B.~Unger.
\newblock Data-driven interpolation of dynamical systems with delay.
\newblock {\em Systems \& Control Letters}, 97:125--131, 2016.

\bibitem{schulze2018data}
P.~Schulze, B.~Unger, C.~Beattie, and S.~Gugercin.
\newblock Data-driven structured realization.
\newblock {\em Linear Algebra and its Applications}, 537:250--286, 2018.

\bibitem{zhang2021factorization}
Q.~Zhang, I.~V. Gosea, and A.~C. Antoulas.
\newblock Factorization of the {L}oewner matrix pencil and its consequences.
\newblock {\em arXiv preprint arXiv:2103.09674}, 2021.

\bibitem{zuyev2007stabilization}
Alexander Zuyev and Oliver Sawodny.
\newblock Stabilization and observability of a rotating Timoshenko beam model.
\newblock {\em Mathematical Problems in Engineering}, 2007.

\end{thebibliography}

\end{document}